\documentclass[11pt]{amsart}
\input{standard.sty}
\usepackage{graphicx,color}

\title{Resolvent estimates with mild trapping}
\author{Jared Wunsch}

\date{\today}

\newcommand{\executeiffilenewer}[3]{%
\ifnum\pdfstrcmp{\pdffilemoddate{#1}}%
{\pdffilemoddate{#2}}>0%
{\immediate\write18{#3}}\fi%
}

\newcommand{\loc}{\text{loc}}
\newcommand{\hamvf}{\mathsf{H}}
\newcommand{\cN}{\mathcal{N}}
\newcommand{\sH}{\mathsf{H}}
\newcommand{\tX}{\widetilde{X}}

\thanks{The author is grateful to Dean Baskin, Hans Christianson,
  Emmanuel Schenck, Andr\'as Vasy, and Maciej Zworski both for
  scientific discussions and for helpful comments on this manuscript.
  The author also thanks David Lannes and the scientific committee of
  the CNRS GDR ``Analyse des Equations aux D\'eriv\'ees Partielles''
  for the opportunity to present the talk on which this exposition is
  based.  This work was partially supported by NSF grant DMS-1001463.}

\begin{document}

\begin{abstract}We discuss recent progress in understanding the effects of
 certain trapping geometries on cut-off resolvent estimates, and thus on the
  qualititative behavior of linear evolution equations.  We focus on
  trapping that is unstable, so that strong resolvent estimates hold
  on the real axis, and large resonance-free regions can be shown to
  exist beyond it.\end{abstract}

\maketitle

\section{Introduction}

Let $(X,g)$ be a Riemannian manifold isometric to $\RR^n$ outside a compact set.  Let $\Lap$ denote the (non-negative) Laplace-Beltrami operator
$$
\Lap=\frac 1{\sqrt{g}} D_i g^{ij} \sqrt{g} D_j,
$$
where $D_i = -\imath \pa_{x_i}.$  The \emph{resolvent}
$$
(\Lap-\lambda^2)^{-1}
$$
is a priori defined as a family of bounded operators $L^2(X)\to
L^2(X),$ as $\lambda$ runs over the upper half-space $\Im \lambda>0.$
When the manifold is precisely $\RR^n$ we can write down the Schwartz
kernel of the resolvent by Fourier methods;
for instance in $\RR^3$ the Schwartz kernel is
$$\frac{e^{i\lambda\smallabs{x-y}}}{4 \pi \smallabs{x-y}}.$$  We thus see
explicitly that this kernel analytically continues to be an
entire function of $\lambda \in \CC.$ It has the modest defect of exponential
growth at infinity, which is remedied if we consider it as a map
$L^2_c(\RR^3)\to L^2_{\loc}(\RR^3).$ The same holds true on $\RR^n$ for all
odd $n.$

Putting geometry back into the picture, we may regard $X$ as a
perturbation of Euclidean space, and an easy application of analytic
Fredholm theory (see, e.g., \cite{Sjostrand-Zworski1} in a much more general
setting) shows that in general, for $n$ odd, $(\Lap-\lambda^2)^{-1}$
continues as a \emph{meromorphic} family of operators $L^2_c(X)\to
L^2_{\loc}(X),$ for all $\lambda \in \CC.$ The poles of this family of
operators are known as the \emph{resonances} of the Laplacian.  We may
equally well think of them as the poles of the \emph{cutoff resolvent
  family}
$$
\chi (\Lap-\lambda^2) \chi: L^2(X)\to L^2(X)
$$
with $\chi \in \CcI(X),$ and we will adopt this point of view in what
follows.  (For a survey of some recent developments in the theory of
resonances complementary to the discussion here, we refer the reader
to \cite{Zworski:Notices}.)

We note that while the new results described in this paper fit (more
or less) in the geometric context described above, much of the history
of the subject is closely tied to the problem where $X$ is replaced by
an exterior domain of a smoothly bounded subset of $\RR^n$ (the
``obstacle problem'') equipped with, say, Dirichlet boundary
conditions.  As the obstacle problem has considerable similarities to that
considered here, we will cite results on the obstacle problem as
motivation without further comment.

We now recall two different estimates in the theory of evolution
equations that follow from an ``appropriate'' understanding of the
resolvent continued to the real axis and beyond.

\subsection{Energy decay for the wave equation}
\begin{figure}
\centering
\def\svgwidth{1in}
\executeiffilenewer{gamma.svg}{gamma.pdf}%
{inkscape -z -D --file=gamma.svg %
--export-pdf=gamma.pdf --export-latex}%
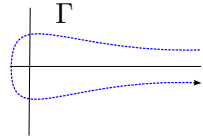%

\label{figure:contour}
\caption{The contour of integration in \eqref{wavecontour}.}
\end{figure}

Let $U(t) ={\sin t \sqrt\Lap}/{\sqrt\Lap}$ be (part of the) wave
propagator.  Then we may write $U(t)$ as a contour integral around the
spectrum of the Laplacian, which is in this situation simply
continuous spectrum along the positive real axis.  Thus, letting
$\Gamma$ denote the contour shown in Figure~\ref{figure:contour}, we
have
\begin{equation}\label{wavecontour}
\begin{aligned}
\chi U(t)\chi &= \frac{1}{2\pi \imath}\int_\Gamma \chi (\mu-\Delta)^{-1}\chi \frac{\sin
t\sqrt{\mu}}{\sqrt\mu}\, d\mu\\
&=-\frac{1}{\pi \imath}\big(\int_{\RR+\imath 0} \chi R(\lambda)\chi
e^{\imath t\lambda} \,
d\lambda+\int_{\RR+\imath 0} \chi R(\lambda)\chi e^{-\imath t\lambda} \, d\lambda\big)
\end{aligned}
\end{equation}

Here we have made the change of variables $\mu=\lambda^2$ and pushed
the contour within $\epsilon$ of the real axis for any $\ep>0$ (hence
the contour is written $\RR+\imath 0$).

Now suppose we are interested in the local large-$t$ asymptotics of a
solution to the wave equation, i.e., of $\chi U(t) f$ for $f$
compactly supported, or equivalently to $\chi U(t) \chi g$ for any $g
\in L^2(X);$ here $\chi$ is again a compactly supported cutoff
function.  For $t \gg 0,$ can can always slide the contour of
integration in the first integral on the last line of
\eqref{wavecontour} upward to $\RR+\imath \nu$ ($\nu>0$) and we get exponential
decay, $O(e^{-\nu t})$ for this term.  The second integral, however,
is trickier, as we would wish to slide this contour \emph{downward}
below the real axis, where we only know that the resolvent is
meromorphic.

In order to deal with this second term, let us now suppose we have:
\begin{enumerate}
\item $\chi (\Lap-(\mu+\imath \nu)^2)^{-1}\chi$ analytic for $-\nu_0 \leq \nu\leq 0.$
\item $\displaystyle\norm{\chi (\Lap-(\mu+\imath \nu)^2)^{-1}\chi}_{L^2}
  \lesssim \ang{\mu}^N$ in this strip.\footnote{Here and henceforth,
    $f\lesssim g$ means $f \leq C g$ for some constant $C.$}
\end{enumerate}

\begin{figure}
\centering
\def\svgwidth{1.5in}
\executeiffilenewer{resfree.svg}{resfree.pdf}%
{inkscape -z -D --file=resfree.svg %
--export-pdf=resfree.pdf --export-latex}%
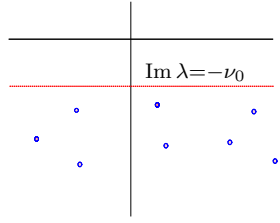%

\caption{A resonance-free strip.}
\end{figure}

If these conditions hold, then we can shift the second contour, at the cost of some powers of $\mu$ or,
equivalently, powers of $\Lap$ and for $K$ sufficiently large get the bound $O(e^{-\nu t})$ for
any $\nu<\nu_0$ for the
smoothed wave operator $\chi \ang{\Lap}^{-K} U(t)\chi.$ Thus we
obtain \emph{exponential decay of local energy} (potentially with derivative loss).

We further remark that in flat space, $N=-1$ in the second condition
and we in fact have no derivative loss in the estimates: whenever
$(D_t^2 -\Lap)u=0$ and $u$ has compactly supported Cauchy data,
$$
\norm{\chi u_t}^2 + \norm{\chi \nabla u}^2\leq C e^{-\nu t} \big(\norm{u_t(0)}^2 + \norm{\nabla u(0)}^2\big)
$$
(where the constant may depend on the size of the support).

\subsection{Local smoothing for the Schr\"odinger equation}

We now consider a second application, which only involves the cut-off
resolvent estimate on the real axis.  This application is to the initial value problem for the time-dependent Schr\"odinger equation:
\begin{equation}\label{scheqn}
\begin{aligned}
\big(\imath^{-1} \frac{\pa}{\pa t}+ \Lap) u&=0,\
u|_{t=0}=u_0.
\end{aligned}
\end{equation}
The estimate on the cut-off resolvent that holds on the real axis in Euclidean space,
\begin{equation}\label{eucres}
\norm{\chi R(\mu)\chi}_{L^2} \lesssim \ang{\mu}^{-1},
\end{equation}
implies the following property of $u:$
$$
\int_0^1 \norm{\chi u}_{H^{1/2}}^2\, dt \lesssim \norm{u_0}_{L^2}^2.
$$
In other words, we have the following mapping property of the propagator:
$$
e^{-\imath t\Lap}: L^2\to L^2([0,1]; H^{1/2}_\loc (X))
$$
(see Constantin-Saut \cite{CS}, Sj\"olin \cite{sjolin}, Vega
\cite{Vega}, Kato-Yajima \cite{Kato-Yajima}, Yajima \cite{Yajima2}).
A proof that \eqref{eucres} gives this estimate proceeds via a standard ``$T T^*$'' estimate which reduces to proving the corresponding estimate on the inhomogeneous equation:
$$
\big(\imath^{-1} \frac{\pa}{\pa t}+ \Lap) v=\chi f,\ u|_{t=0}=0 \Longrightarrow \norm{\chi u}_{L^2_t H^{1/2}_x} \lesssim \norm{\chi f}_{L^2_t H^{-1/2}_x}.
$$
This inhomogeneous estimate, in turn, follows from Fourier transform and \eqref{eucres}.

\subsection{Semi-classical rescaling, and estimates in non-trapping geometries}

Having established that cut-off resolvent estimates on and below the
real axis are of key importance in the study of the wave and
Schr\"odinger equations, we now pursue the question of when they hold
in non-Euclidean geometries.  To begin, we note that both the validity
of the estimate \eqref{eucres} and also of the existence of a
pole-free strip with polynomial resolvent bounds are
\emph{high-frequency} questions: they manifestly hold in any compact
range of $\Re \lambda,$ hence they hinge only on the
asymptotics of the resolvent as $\Re \lambda \to +\infty.$
Consequently, we now introduce a \emph{semi-classical rescaling} of
the problem: we set
$$
\lambda^2=\frac{z}{h^2}
$$
and study the operator family
$$
P_h(z) \equiv h^2 \Lap-z
$$
and its inverse
$$
R_h (z) \equiv (h^2 \Lap-z)^{-1}.
$$
We now note that
the existence of a resonance
free strip in $\lambda$ is equivalent to pole-free region in $z$
for $\chi R_h(z)\chi$ of form
$$
[1-\delta,1+\delta] + i [-\nu_0 h, 0],
$$
and the desired polynomial growth estimate in $\Re \lambda$ becomes a polynomial estimate in $h^{-1}.$
Likewise the free resolvent estimate on the real axis \eqref{eucres} rescales to:
$$
\norm{\chi R_h(z) \chi} \lesssim h^{-1},\quad z \in [1-\delta,1+\delta].
$$

In order to explain a broad setting in which these estimates
generalize, we now introduce the notion of \emph{trapped set.}
Let $\hamvf_p$ denote the Hamilton vector field of
$p=\sigma_h(P),$ i.e., the generator of geodesic flow in the cosphere
bundle. Then the trapped set $K$ is defined by
$$
\rho \in K \Longleftrightarrow \exp(t\hamvf_p)
 \text{ remains in a compact set for all }  t.
$$
We say the metric is \emph{non-trapping} if $K=\emptyset.$

The following is the classical result on estimates for non-trapping metrics:
\begin{proposition}
If the metric is nontrapping then there is a resonance free region of the form
$$
[1-\delta,1+\delta] +\imath [-\nu_0 h\smallabs{\log h}, 0]
$$
and the free
resolvent estimate holds: $$\norm{\chi R_h(z) \chi} \lesssim h^{-1},\ z \in [1-\delta,1+\delta].$$
\end{proposition}
Indeed, a somewhat stronger statement holds on the distribution of
resonances: for $h$ small, the region $[1-\delta,1+\delta] +\imath [-\nu h\smallabs{\log h}, 0]$
contains only finitely many resonances, for \emph{any} value of $\nu.$

These estimates have a lengthy history in both the obstacle case and
the metric setting considered here.  The deduction of exponential
energy decay from the eventual smoothness of solutions to the wave
equation with compactly supported data (a very weak form of Huygens'
Principle) follows from work of Lax-Phillips \cite{Lax-Phillips1} and Vainberg
\cite{Vainberg:Exterior,Vainberg:Asymptotic}.  (For further references
on the propagation of singularities arguments leading to these
estimates in the obstacle case, we refer the reader to
\cite{Tang-Zworski1}; in the metric case dealt with here, it simply
follows from H\"ormander's theorem on propagation of singularities for
operators of real principal type \cite{MR48:9458}.)

The estimate on the real axis is sharp: Ralston
\cite{Ralston:Localized} showed (in the obstacle setting) that the
estimate on the real axis can \emph{only} hold if the trapped set is
empty.  Indeed, if there is \emph{stable} trapping, e.g.\ a closed
elliptic geodesic, then it is known that the resolvent grows
on $\RR,$ with resonances approaching the real axis, as $\Re \lambda \to
 +\infty$ (see e.g., \cite{Stefanov-Vodev1,Stefanov-Vodev2}, \cite{Tang-Zworski2},
\cite{Stefanov:Quasimodes}). The minimal size
of the resonance free strip was established by Burq \cite{Burq1}, who
showed that there must always be a region of the form $\Im
\lambda>-Ce^{-\ep\smallabs{\lambda}}$ free of poles of
$(\Lap-\lambda^2)$ and consequently there is energy decay for the wave
equation at the rate $\log(2+t)^{-k}$ if we accept loss of $k$
derivatives.

\section{Resolvent estimates in the presence of mild trapping}

\subsection{Normally hyperbolic trapping}
What, then, if there are trapped orbits, i.e., $K \neq \emptyset,$ but
the trapping is at least \emph{unstable}?  The classical example of Ikawa sheds
considerable light on this problem: in the case of the exterior
problem for two strictly convex obstacles, there is a unique trapped
orbit consisting of the single orbit bouncing back and forth between
them; this orbit is highly unstable, and Ikawa was able to show that
in this case there is a strip $\Im z \geq -\nu_0 h$ that is free of
resonances, and indeed to derive asymptotics for the locations of the
resonances.  A refinement of Tang-Zworski \cite{Tang-Zworski1} and Burq
\cite{MR2066943} shows that while the non-trapping estimate cannot
hold on the real axis, we still have:
$$
R_h(\lambda) \lesssim \frac{\smallabs{\log h}}{h}, \quad \lambda \in [1-\delta,1+\delta].
$$
In other words, we lose only a factor of $\abs{\log h}$ relative to
the non-trapping estimate.  This implies, among other things, that the
local smoothing estimate for the Schr\"odinger propagator only barely
fails: for any $\ep>0,$
$$
e^{-\imath t\Lap}: L^2 \to L^2([0,1]; H_{\loc}^{1/2-\ep}).
$$
In the boundary-less case, the analogous situation in
which we have a single unstable hyperbolic orbit (e.g., the unique closed
geodesic on the surface of revolution shown in
Figure~\ref{figure:revolution} below)
was analyzed by
Christianson \cite{Christianson1}.  More generally, the author and
Zworski \cite{Wunsch-Zworski:normhyp} have recently studied the situation of \emph{normally
  hyperbolic trapping}, in which there is a whole manifold of unstable
trapped orbits.  We now describe these results.

One motivation for considering this geometric set-up comes from the
\emph{Kerr black hole.}  This is a family of Lorentzian metrics which
solve the Einstein equations and describe rotating black holes.  The
trapped set of null-geodesics (known as the ``photon sphere'' in the
special, spherically symmetric, case of the Schwarzschild metric) has the dynamical structure
described below, and thus the results of \cite{Wunsch-Zworski:normhyp}
address the obstructions to decay posed by the trapped rays.
Considerable geometric and analytic difficulties are present in the
study of wave decay on Kerr backgrounds that are not addressed by
\cite{Wunsch-Zworski:normhyp} however; see, e.g.
\cite{Tataru2},
\cite{Tataru-Tohaneanu1} \cite{DSS1}, \cite{Dafermos-Rodnianski3},
\cite{Andersson-Blue1}, \cite{FSKY,FSKY:erratum} for results along these
lines and for further references.
  In the
related setting of Kerr-de-Sitter black hole metrics\footnote{These
  metrics describe rotating black holes in a universe with positive
  cosmological constant.} and their perturbations, some of the
obstructions to decay posed by the asymptotically Euclidean end of the
de Sitter space are not present.  In the case of
slowly-rotating\footnote{Slowly rotating means that the angular
  momentum parameter in the metric, usually denoted $a,$ is taken to
  be small.}
Kerr-de-Sitter metrics,
Dyatlov \cite{Dyatlov:KdS1,Dyatlov:KdS2} used the results from
\cite{Wunsch-Zworski:normhyp} presented
below to show that wave equation solutions decay exponentially to a
constant value; Vasy \cite{Vasy:KdS} was able to extend such
decay estimates to a wide class of perturbations of Kerr-de-Sitter
space, again using the results of \cite{Wunsch-Zworski:normhyp} to cope
with the trapped set.

The dynamical hypotheses employed in \cite{Wunsch-Zworski:normhyp} are as follows.
For simplicity of exposition here we stick to the geometric set-up
described above, where we consider the operator $h^2\Lap$ on the
manifold $(X,g)$
with $g$
Euclidean outside a compact set.
Let $\varphi^t$ denote the corresponding geodesic flow.
Let
$r$ denote the distance function to a fixed point in $X$ and locally
define the backward/forward trapped sets by:
$$
\Gamma_{\pm}=\{\rho: \lim_{t \to \mp\infty} r(\varphi^t(\rho)) \neq \infty\}.
$$
The the \emph{trapped set} is simply
$$
K=\Gamma_+\cap \Gamma_-.
$$

We make dynamical assumptions as follows:
\begin{enumerate}
\item \label{gammas}$\Gamma_\pm$ are codimension-one 
smooth
manifolds intersecting transversely at $K.$ (It
is not difficult to verify that $\Gamma_\pm$ must then be coisotropic and
$K$ symplectic.)

\item \label{normflow} The flow is hyperbolic in the normal directions
to $K$ within the energy surface: letting the superscript $1$
denote the intersection of the set in question with the unit cosphere
bundle (e.g., $K^1=K\cap S^*X$) then there exist subbundles $E^\pm$ of
$T_{K^1} \Gamma_\pm^1$ (the ``unstable/stable subspaces'') such that 
$$
T_{K^1} \Gamma_\pm^1=T K^1 \oplus E^\pm,
$$
where $$d\varphi^t:E^\pm \to E^\pm$$ and
there exists $\theta>0$ such that on $K^1,$
\begin{equation}
\label{eq:normh}
\| {d\varphi^t(v)} \|  \leq C e^{-\theta | t | } \| {v} \| \
\text{for all } v \in E^{\mp},\ \pm t \geq 0.
\end{equation}
\end{enumerate}

These assumptions can be verified directly
(see \cite{Wunsch-Zworski:normhyp}) for the trapped set of a
slowly rotating Kerr black hole\footnote{See also \cite{Vasy:KdS} for
  analogous discussion of the Kerr-de-Sitter metric.} (i.e.\ when the
angular momentum parameter $a$ is small---see
\cite{Wunsch-Zworski:normhyp} for the explicit form of the metric) but
they are not stable under perturbations, hence do not obviously apply
to perturbations of Kerr metrics.  However, the bicharacteristic flow
for the Kerr metric in fact satisfies a more
stringent (and well-studied) hypothesis that \emph{is} stable under
perturbation, and that implies the dynamical hypotheses above.  In
particular, the standard dynamical notion of \emph{$r$-normal
  hyperbolicity} implies items \eqref{gammas} and \eqref{normflow},
and \emph{is} stable under perturbations, modulo possible loss of
derivatives.  This hypothesis says, roughly speaking, that along $K$
the tangent space of the normal bundle to $K$ splits into subbundles
$E^\pm$ which are exponentially expanded/contracted by the flow, while
the expansion/contraction of directions along $TK$ is much milder by
comparison (e.g., polynomial in the case of the Kerr metric).  It is
then a profound theorem of Fenichel \cite{Fenichel:Persistence} and
Hirsch-Pugh-Shub \cite{Hirsch-Pugh-Shub1} that there exist
stable/unstable manifolds $\Gamma_\pm$ tangent to $E^\pm$ and
satisfying the dynamical hypotheses above, and moreover that these
stronger hypotheses are structurally stable.\footnote{Albeit in
  general, perturbing results in loss of differentiability of the
  manifolds $\Gamma_\pm$ and $K;$ however for any desired
  \emph{finite} degree of differentiability, shrinking the
  perturbation preserves that degree of smoothness, which suffices for
  our purposes.}

Our main result on the resolvent for normally hyperbolic trapping is
as follows:
\begin{theorem}{(Wunsch-Zworski \cite{Wunsch-Zworski:normhyp})}\label{theorem:normhyp}
Under the normal hyperbolicity hypotheses described above, there
exist $\delta, \nu_0>0$ such that $(h^2\Lap-z)^{-1}$ is
meromorphic in the region
$$
[1-\delta,1+\delta] +\imath [-\nu_0 h, 0]
$$
and the following polynomial resolvent estimates on the hold: on the
real axis, we have $$\norm{\chi R_h(\lambda) \chi} \lesssim
\frac{\smallabs{\log h}} h,\ \lambda \in [1-\delta,1+\delta],$$ while
  below the axis there exists $k>0$ such that
$$\norm{\chi R_h(z) \chi} \lesssim
h^{-k},\ z\in [1-\delta,1+\delta] +\imath [-\nu_0 h, 0].$$
\end{theorem}

In the case when the operator has real analytic coefficients 
a more precise resonance free region was obtained by 
G\'erard-Sj\"ostrand \cite{Gerard-Sjostrand:Lyapunov}, albeit without the polynomial bound
on the resolvent that makes contour deformation (with loss of
derivatives) feasible to obtain exponential energy decay for solutions
to the wave equation.  A similar estimate on the optimal width of the
resonance free strip in the setting considered here, and with much reduced
smoothness assumptions on the stable/unstable manifolds $\Gamma_\pm,$ has been
recently announced by Nonnenmacher-Zworski \cite{Zworski:Private}.

Estimates like Theorem~\ref{theorem:normhyp} have also been proved in
the case when $K$ is \emph{fractal} rather than a smooth manifold:
provided the set is sufficiently filamentary, in a manner measured by
a \emph{topological pressure} condition, an estimate of the same form
was previously proved by Nonnenmacher-Zworski
\cite{Nonnenmacher-Zworski1}; in \cite{N-Z:Semiclassical} these
authors also obtain more precise resolvent bounds below the real axis.
(See also Christianson \cite{Christianson:Cutoff} and Datchev \cite{Datchev:Local} for
applications to Schr\"odinger local smoothing in the presence of
various assumptions on the ends of $X$.)

\subsection{Degenerate hyperbolic trapping}

Since stable trapping, for instance by elliptic closed geodesics, generates
exponentially growing resolvent estimates, while unstable trapping
with hyperbolic dynamics 
leads only to a loss of $\smallabs{\log h}$ on the real axis, the
reader might wonder whether there are situations in which a nontrivial
polynomial loss occurs.  This phenomenon was demonstrated by the
author and Hans Christianson \cite{Christianson-Wunsch1} in the following
situation:
We consider the manifold $X = \RR_x \times \RR_\theta / 2 \pi
\ZZ$, equipped with a metric of the form
\[
ds^2 = d x^2 + (1 + x^{2m})^{1/m} d \theta^2,
\]
Here the trapped set is again the closed orbit where $x=0.$  If $m=1,$
this is again the classic case of a single trapped hyperbolic orbit,
but if $m$ is an integer greater than $1$ the hyperbolicity
\emph{degenerates} with the stable and unstable manifold no longer
intersecting transversely.  Christianson and the author showed the
following:
\begin{theorem}{(Christianson-Wunsch \cite{Christianson-Wunsch1})}
For $m$ an integer at least $2,$
$$
\norm{\chi R_h(\lambda)\chi}\lesssim h^{-2m/(m+1)},\quad \lambda \in [1-\delta,1+\delta]
$$
and this estimate is sharp.
\end{theorem}
Consequently, the optimal local smoothing estimate in this situation entails a
nontrivial loss of derivatives:
$$
e^{-\imath t\Lap}: L^2 \to L^2([0,1]; H_{\loc}^{1/(m+1)}),
$$
and we may obtain a polynomial rather than exponential rate of energy decay for the wave
equation.

\subsection{Trapping by cone points}

While metrics with trapped rays always have losses in resolvent
estimates on the real axis and in Schr\"odinger local smoothing (see
Doi \cite{Doi1} and Burq-Bony-Ramond \cite{B-B-R} as well as Ralston
\cite{Ralston:Localized} in the obstacle case), we now present a
situation involving wave propagation on singular manifolds in which a
weaker form of ``diffractive'' trapping results in (almost) no loss whatsoever.
This is the setting of \emph{manifolds with cone points.}

In joint work with Dean Baskin, the author considers a manifold $X$ of
dimension $n$ with conic singularities that is isometric to $\RR^n$
outside a compact set.  Here a manifold with conic singularities means
a manifold $X$ with compact boundary $Y=\pa X$ equipped with a metric
that can be brought to the form
$$
g=dx^2 + x^2 h
$$
in a neighborhood of $Y,$ with $h$ a smooth tensor restricting to
give a metric on $Y.$  The metric thus degenerates at the
boundary, and each boundary component becomes a cone
point.

In considering trapping, we must now distinguish two kinds of geodesic
on $X.$  Following Melrose-Wunsch \cite{Melrose-Wunsch1},  we let a
\emph{diffractive geodesic} denote a geodesic that can enter and leave
a cone point (i.e., a single boundary component) along \emph{any} pair
of geodesics.  By contrast, a \emph{geometric geodesic} is one that is
restricted to enter and leave a single cone point along a pair of
geodesics that are connected by a
geodesic of length $\pi$ inside $Y$ (with respect to the metric $h\rvert_{Y}$).  In \cite{Melrose-Wunsch1},
the geometric geodesics are shown to be exact those that are (locally)
approximable by ordinary geodesics in $X^\circ.$

We impose geometric hypotheses
that \begin{enumerate}
\item The flow
along geometric geodesics is non-trapping.
\item 
No three cone points are collinear along a geometric geodesic.
\item
No two cone points are conjugate to one another.
\end{enumerate}
(The last condition can be most easily interpreted as saying that the geodesic
flowouts from two different cone points intersect transversely.)

Then the following estimate on the resolvent holds:
\begin{theorem}{(Baskin-Wunsch \cite{BW1})}\label{theorem:conedecay}
For $\chi \in \CcI(X),$ there exist $\delta, \nu_0>0$ such that the cut-off resolvent
$$
\chi (h^2\Lap-z)^{-1}\chi
$$
can be analytically continued from $\Im z >0$ to the region
$$
z\in [1-\delta,1+\delta]+\imath [-\nu_0 h \smallabs{\log h}, 0]
$$
and for some $C,T>0$ enjoys the estimate
$$
\norm{\chi (h^2\Lap-z)^{-1}\chi}_{L^2\to L^2} \leq \frac{C}{h} e^{T \smallabs{\Im z}/h}
$$
in this region.
\end{theorem}
This implies that the trapping of waves induced by diffraction among
cone points has only a very weak effect on the decay of solutions to
the wave equation on such a manifold: in particular, solutions to the
wave equation on $\RR\times X$ enjoy exponential local energy decay.
Likewise, the resolvent estimate on the real axis (which is the same
as would be obtained with \emph{no} trapped rays) implies that the
Schr\"odinger propagator $e^{-\imath t\Lap}$ enjoys the local smoothing
estimate with no loss, mapping $L^2\to L^2([0,1]_t; H^{1/2}_{\loc}(X)).$

The theorem is proved by first establishing a result on weak escape of
singularities for solutions to the wave equation.  We recall that in
odd dimensional Euclidean space, Huygens' Principle tells us that any
solution to the wave equation with compactly supported Cauchy data is
eventually zero in any fixed compact set $F$; in a non-trapping
perturbation of Euclidean space (of any dimension), the solution does
not in general vanish after long time but instead, by propagation of
singularities, is eventually in $\CI(F).$ In the conic situation, by
contrast, we expect that singularities may persist forever in a
compact set, diffracting back and forth among the cone points, but it
turns out that they become milder and milder as time passes:
\begin{theorem}{(Baskin-Wunsch \cite{BW1})}\label{theorem:conesmoothing}
For any
$r$ there exists $T_r$ such that for $\abs{t}>T_r,$ for all $s,$
$$
\chi \cos t\sqrt{\Lap}\chi: \mathcal{D}_s\to \mathcal{D}_{s+r}.
$$
\end{theorem}
\noindent (Here $\mathcal{D}_s$ is the domain of $\Lap^{s/2}$ and agrees with
$H^s$ away from the cone points.)  One can then deduce
Theorem~\ref{theorem:conedecay} from
Theorem~\ref{theorem:conesmoothing} by adapting an argument of
Vainberg \cite{Vainberg:Asymptotic} that employs the Fourier transform
of the outgoing singular parts of the solution as an initial
parametrix for the resolvent, solving away the resulting errors by
using the free wave operator (for which the resolvent estimate
assuredly holds).

One especially interesting and elementary setting in which our result
on resonance-free regions applies is that of exterior domains to
polygons in $\RR^2.$  We can \emph{double} such domains by gluing
together two copies along their boundary to obtain a manifold with
cone points; solving the wave equation with Dirichlet/Neumann boundary
conditions is equivalent to working with odd/even solutions on the
doubled manifold.  Provided the domain is
nontrapping in the sense that billiard trajectories missing the cone
points escape to infinity, we can thus show that there is a
logarithmic region free of resonances for the obstacle problem with
either boundary condition.\footnote{There is some mild technical
  complication from the fact that the manifold has two Euclidean ends,
  rather than just one as in our hypotheses, but the same strategy of
  proof applies.}

We remark that the size of the resonance free region in
Theorem~\ref{theorem:conesmoothing} is probably sharp: in the related
situation of the exterior problem in $\RR^2$ outside two strictly
convex analytic obstacles, one having a corner, Burq \cite{Burq:Coin}
has shown that resonance poles for $(\Lap-\lambda^2)^{-1}$ lie
(asymptotically) along curves of the form $\Im \lambda=C\abs{\Re \lambda} +D.$
The poles are thus at distance $O(h
\smallabs{\log h})$ from the real axis after semiclassical rescaling.  (We also refer the reader to
Zworski \cite{Zworski4} where analogous logarithmic strings of poles are
shown to arise from finite order singularities of a one-dimensional
potential, thus substantiating heuristics from
Regge \cite{Regge:Analytic}.)  Consequently, the only loss in this setting of
``diffractive trapping'' relative to the non-trapping case is that we
expect that some region of the form $z\in [1-\delta,1+\delta]+\imath [-\nu_0 h \smallabs{\log h}, 0]$
probably contains infinitely many resonances in the former case, while
this is not so in the latter.

\section{From resolvent estimates to damped waves}

We now turn to a problem related to, but crucially distinct from, that
of resolvent estimates and decay for solutions to the wave equation:
we now consider the \emph{damped wave equation} on a \emph{compact,} connected
Riemannian manifold $X.$  We fix a nonnegative ``damping function''
$a(x) \in \CI$ and consider the initial value problem:\footnote{Recall
  that the Laplacian employed here is the \emph{positive} operator $\Lap=d^*d.$}
\begin{equation}
\label{eq:dwe}
\bigg\{ \begin{array}{l}
\left( \partial_t^2 + \Delta_g + 
a(x) \partial_t \right) u(x,t) = 0, \\
u(x,0) = u_0\in H^1(X), \quad \partial_t u(x,0) = u_1\in H^0(X).
\end{array}
\end{equation}
It is well known \cite{MR94b:93067}, \cite{MR97i:58173} that if $a>0$
somewhere, then any solution tends to zero as $t\to +\infty.$  Moreover
Lebeau has shown (see \cite{MR97i:58173}) that the energy of the waves,
$$
E(u,t)=\frac{1}{2}\int_X \big( | \partial_t u |^2 + | \nabla u |^2
\big) dx 
$$
decays at a logarithmic rate provided the data lie in a space with
higher regularity than the energy space; this is of course analogous
to the results discussed above for energy decay of solutions to the wave equation
on noncompact manifolds in the presence of trapping.

To understand decay rates for the damped wave equation in further
detail, we make a semiclassical reduction: Fourier transforming, it
turns out that we can replace the damped wave operator by an operator
of the form
\begin{equation}\label{damped}
h^2\Lap + \imath h \sqrt{z} a-z;
\end{equation}
exponential energy decay then
follows from the existence of a pole-free region of the form
\begin{equation}\label{polefreeregion} z\in [1-\delta,1+\delta]+\imath
  [-\nu_0 h,0]\end{equation} together with polynomial
bounds for the inverse of this non-self-adjoint operator; smaller
regions of analyticity yield sub-exponential decay rates.  In
the results discussed here, the $z$ dependence of the middle term in \eqref{damped} will
be treatable as a perturbation in the operator above, and we study the family
\begin{equation}\label{dampingresolvent}
(h^2\Lap + \imath h a-z)^{-1}
\end{equation}
which is at least formally quite similar to the resolvent studied in
previous sections.

Now in analogy with the trapped set in noncompact problems, we introduce the \emph{undamped set}
\begin{equation}
\cN=\{\rho\in S^*X: \forall t\in\RR, a\circ e^{t\sH_p}(\rho)=0 \}
\end{equation}
consisting of points that do not encounter the interior of the damping region under
the forward- or backward-bicharacteristic flow.  If $\cN=\emptyset$ we
are said to be in the situation of \emph{geometric control.}
Classic results of Rauch-Taylor \cite{RaTa75}, Bardos-Lebeau-Rauch
\cite{MR94b:93067} and Lebeau \cite{MR97i:58173} show that in the setting of geometric control, we
do have a pole-free region of the form \eqref{polefreeregion} and
exponential decay of energy.  By contrast, if $\cN\neq
\emptyset,$ it is known that exponential decay cannot hold without
derivative losses.  The situation is thus closely analogous to that
described above for resolvent estimates, leading one to wonder if
there is a link between the resolvent estimates for a noncompact
manifold $\tX$ with trapped set $K$ and estimates for
\eqref{dampingresolvent} on a compact manifold $X$ with undamped set $\cN
\equiv K,$ as long as the two manifolds are isometric in a
neighborhood of the undamped/trapped set.  In joint work with Hans
Christianson, Emmanuel Schenck, and Andr\'as Vasy, the author has
proved that this is so: provided one has polynomial resolvent bounds
on the real axis,
one can ``glue'' these bounds to obtain the exact same estimate for
the damped resolvent family \eqref{dampingresolvent}.  While our
motivation comes from gluing on non-compact ends, the following
theorem is in fact phrased in terms of the resolvent for the problem with
``complex absorbing potential'' $\imath a$ replacing the damping
coefficient $\imath h a;$ such a potential has the effect of
annihilating semiclassical wavefront set along (forward)
bicharacteristics passing through it.  Estimates with non-compact ends
and with complex absorbing potential are known to be
equivalent in a wide variety of geometric settings by a recent gluing
theorem of Datchev-Vasy \cite{Datchev-Vasy:Gluing}.  Phrasing the
resolvent estimate in this way allows us to stay within the setting of
compact manifolds, however:
\begin{theorem}{(Christianson-Schenck-Vasy-Wunsch \cite{CSVW})}
\label{theorem:damping}
Assume  that for some $\delta \in (0,1)$
fixed and $k \in \ZZ$, there is a function $1 \leq \alpha(h) = O( h^{-k})$ such that
\[
\| (h^2\Delta_g +\imath a-z)^{-1} \|_{L^2 \to L^2} \leq \frac{\alpha(h)}{h},
\]
for $z \in [1-\delta, 1+
\delta]$.  Then there exist $C, \nu_0>0$ such that 
\[
\| (h^2\Delta_g+\imath ha -z )^{-1} \|_{L^2 \to L^2} \leq C \frac{\alpha(h)}{h},
\]
for $z \in [1-\delta,1+\delta] + \imath [-\nu_0 h/\alpha(h) , 0 ]$.
\end{theorem}

\begin{figure}
\centering
\def\svgwidth{2in}
\executeiffilenewer{dampedpeanut.svg}{dampedpeanut.pdf}%
{inkscape -z -D --file=dampedpeanut.svg %
--export-pdf=dampedpeanut.pdf --export-latex}%
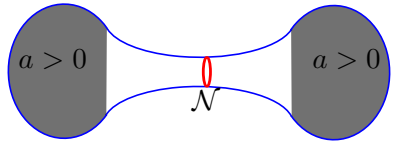%

\caption{The compact manifold $X$ with damping.}
\end{figure}
\begin{figure}
\centering
\def\svgwidth{1.9in}
\executeiffilenewer{noncompact.svg}{noncompact.pdf}%
{inkscape -z -D --file=noncompact.svg %
--export-pdf=noncompact.pdf --export-latex}%
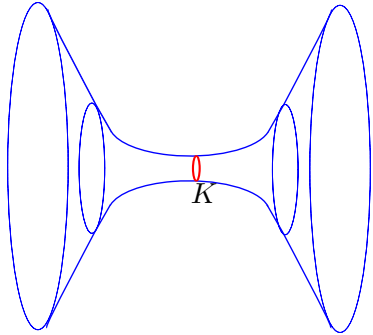%

\caption{The noncompact manifold $\tX.$}
\label{figure:revolution}
\end{figure}
An essential new ingredient in the proof is a recently obtained
estimate of Datchev-Vasy \cite{Datchev-Vasy:Trapped} which improves
the microlocal resolvent estimates available at points which are
backwards- or forwards-trapped relative to those which are trapped
along both directions of the flow.

Theorem~\ref{theorem:damping} has applications in all the examples
discussed above (normally hyperbolic
trapping, fractal hyperbolic trapping, degenerate hyperbolic trapping)
giving various subexponential decay rates for solutions to the damped
wave equation in these settings; we refer the reader to \cite{CSVW}
for further details, including the explicit decay rates thus
obtained.  It is at least sometimes \emph{sharp}, e.g.\ in the case of
a surface of rotation with $\cN$ given by an single undamped
hyperbolic orbit, according to work of Burq-Christianson
\cite{Burq-Christianson:Peanut}.

\bibliography{all}{}
\bibliographystyle{plain}
\end{document}